\documentclass[12pt]{article}

\usepackage{amsmath, amsfonts, amssymb, amsthm} \parskip1pt

\def\cO{ {\cal O}}

\def\bC{{\mathbb  C}}

\def\bN{{\mathbb N}}

\def\bR{{\mathbb  R}}
\def\bZ{{\mathbb  Z}}

\def\bX{{\bf X}}
\def\bp{{\bf p}}
\def\bq{{\bf q}}
\def\by{{\bf y}}
\def\bY{{\bf Y}}

\def\bZ{{\mathbb Z}}

  \def\pro{\noindent {\bf{Proof. }}}

\def\card{{\rm card }}
 \def\Card{{\rm Card }}

\def\build#1_#2^#3{\mathrel{\mathop{\kern 0pt#1}\limits_{#2}^{#3}}}

\def\smallsquare{\vbox{\hrule\hbox{\vrule height 1 ex\kern 1
ex\vrule}\hrule}}
\def\cqfd{\hfill \smallsquare\vskip 3mm}

\def\cB{{\cal B}}

\def\cE{{\cal E}}

\def\mat{{\rm Mat}_{n,m}(\bR)}
\def\Mat{{\rm Mat}_{n,m+n}(\bR)}
\def\Matmax{{\rm Mat}_{n,m+n}^{{\rm rank }\,\, n}(\bR)}
\def\SLZ{{\rm SL}_{m+n}(\bZ)}
\def\SL{{\rm SL}_{2}(\bZ)}

\newtheorem{thm}{Theorem}

\newtheorem*{prop}{Proposition}
\newtheorem{prob}{Problem}

\newtheorem*{Kronecker}{Theorem {\rm (Kronecker,  Dani-Raghavan)}}
\newtheorem*{Erdos}{Theorem {\rm (Erd\H{o}s)}}

\date{ }
\begin{document}

\title{ On Kronecker's  density theorem, primitive points and orbits of matrices }
\author{Michel Laurent \;\; }
\maketitle

{\small
\noindent{\bf Abstract:} 
We discuss recent quantitative results in connexion with  Kronecker's  theorem on the density of subgroups in $\bR^n$
and with  Dani and Raghavan's theorem on the density of orbits in the spaces of frames. We also propose several related problems. The  case of the natural linear action of the unimodular group $\SL$ on the real plane is investigated more closely.
We then establish an intriguing link between the configuration of (discrete) orbits of primitive points and the rate of density of dense orbits.
 }

\vskip 7mm

\section{ Introduction}

Let $m$ and $n$ be positive integers and let $\Theta\in \mat$ be an $n \times m$ matrix with real entries.
We associate to $\Theta$ the subgroup
$$
\Lambda = \Theta \bZ^m+ \bZ^n \subset \bR^n
$$
generated over $\bZ$ by the $m$ columns of $\Theta$ and by $\bZ^n$ and its subset of  primitive points 
$$
\Lambda_{\rm prim} = \Big\{ \Theta \bq + \bp \, ; \, \bq \in \bZ^m , \bp\in \bZ^n\quad {\rm with }\quad \gcd(\bq,\bp)=1\Big\}.
$$
 Put also 
 $$
 X=(\Theta,I_n) \in \Mat ,
 $$
  where $I_n$ is the identity matrix in ${\rm Mat}_{n,n}(\bR)$, so that we can write
$$
\Lambda = X \bZ^{m+n} \quad {\rm and}\quad \Lambda_{\rm prim} = X P(\bZ^{m+n}).
$$
Here and throughout the article,  $P(\bZ^{d})$ denotes the set of  primitive points in $\bZ^{d}$, that is the set of integer $d$-tuples with coprime coordinates.

Then we can state the following  criterion  of density:  
 \begin{Kronecker}
The following assertions  are equivalent.
\\
{\rm(i)}  The $n$ rows of the matrix $\Theta$ are $\bZ$-linearly independent vectors in $\bR^m/\bZ^m$.
\\
{\rm(ii)}    The subgroup $\Lambda$ is dense in $\bR^n$.
\\
{\rm(iii)}  The set $\Lambda_{\rm prim}$ is dense in $\bR^n$.
\\
{\rm(iv)}   The   orbit $X\SLZ $ of $X$ under  the action of the unimodular group $\SLZ$ by right matrix multiplication  is dense in $\Mat$.
\end{Kronecker}

The equivalence of properties (i) and (ii) is the content of the classical  Kronecker density theorem, see for instance \cite{Cas}, 
 while the equivalence of (i) and (iv) is due to Dani and Raghavan,  
see Theorem 3.4 in \cite{DaRa}. 
The implications (iv) $\Rightarrow$(iii) $\Rightarrow$ (ii) are clear;   observe  that upon identifying an $n\times (m+n)$ matrix with the $(m+n)$-tuple of its columns, 
we have the inclusions
$$
X\SLZ  \subseteq \left(X P(\bZ^{m+n})\right)^{m+n}=(\Lambda_{\rm prim})^{m+n}\subset (\bR^n)^{m+n} \simeq \Mat ,
$$
since the    column vectors  of any  matrix in $\SLZ$ belong to $P(\bZ^{m+n})$.

Quantitative  results of density for the subgroup $\Lambda$ are well understood. Assuming that the rows of $\Theta$ satisfy some measure of 
$\bZ$-linear independence modulo $\bZ^m$,
we can control by a transference principle the quality of the approximation of any point in $\bR^n$ by elements of $\Lambda$. We refer to Chapters III and V of Cassels'
monograph \cite{Cas} for this classical issue of Diophantine approximation, and to \cite{BuLa} for a formulation in terms of exponents of approximation.
The inhomogeneous variants  \cite{ScA,Sp} of the metrical Khintchine-Groshev theorem may as well be considered as  quantitative versions of the assertion (ii)
for a generic matrix $\Theta$.

The purpose of this  article is to exhibit  some  effective density results in connexion with the assertions (iii) and (iv), both for a given matrix $\Theta$ 
and for $\Theta$  generic. 
Our knowledge concerning quantitative versions of (iii) and (iv) is much  more limited and 
 another goal of the paper is to formulate some related questions, or conjectures, which may lead to further improvements on these two issues.  
We display in Section 2 a metrical theory involving primitive points, in a generalized meaning, obtaining results refining the
 Khintchine-Groshev theorem. Regarding (iv), even the generic exponent of Diophantine approximation (the analogue of Dirichlet's exponent $m/n$ for $\Lambda$) remains unknown for any value of the dimensions $m$ and $n$. We propose in Section 3 a conjectural value for this critical exponent which is motivated by
recent works \cite{GGNA,GGNB} due to Ghosh, Gorodnick and Nevo in a much more general framework. Sections 4 and 5 are devoted to the case $m=n=1$. We 
 develop a conditional  approach to reach the expected exponent $1/2$ based on  some hypothesis concerning  the repartition in the plane of  truncated
  $\SL$-orbits  of  primitive integer points,  a problem which may have  an independent interest.

\section{Metrical theory for  $\Lambda_{\rm prim}$}

We have obtained in \cite{DLN} metrical statements in the style of the Khintchine-Groshev theorem involving primitive points 
in a refined sense. In this section, we describe and discuss our results.

The set-up is as follows. Let $\Pi$ be a partition of the set 
$
\{1, \dots, m+n\} = \coprod_{j=1}^t \pi_j.
$
 We assume that all  components  $\pi_j$ of $\Pi$ have cardinality at least $n+1$. We then define $P(\Pi)$
to be the set of integer $(m+n)$-tuples $(v_1, \dots , v_{m+n})$ such that 
$\gcd(v_i)_{i\in \pi_j} = 1$  for all $ j=1, \dots, t$.
Note that the trivial partition  (the one having only  one component)  satisfies the assumption  and  that we have $P(\Pi) = P(\bZ^{m+n})$ in that case. 
Moreover $P(\Pi)\subseteq P(\bZ^{m+n})$ for any relevant partition $\Pi$. Now let $\psi : \bN \mapsto \bR^+$ be a  function with positive values. 
We assume that the map $x \mapsto x^{m-1}\psi(x)^n$ is non-increasing and that
the series $\sum_{ \ell\ge 1} \ell^{m-1}\psi(\ell)^n$ diverges; the latter assumption  being the necessary and sufficient condition occurring in the 
 statement of the classical Khintchine-Groshev theorem. 
 
 In the following two theorems, we implicitely assume that the above hypotheses on
 $\Pi$ and $\psi$ are satisfied. The  symbol $| \cdot |$  indicates the supremum norm in $\bR^n$ and in $\bR^m$.
 Our first result is a doubly metrical statement.
 
 \begin{thm}
For almost every pair $(\Theta,\by)\in \mat\times\bR^n$, there exist 
infinitely many points $(\bq,\bp)\in P(\Pi) $ 
such that
$$
| \Theta \bq + \bp - \by |\le \psi(| \bq |). 
$$
 \end{thm}
  \begin{prob}
Fix arbitrarily the target point $\by \in \bR^n$. Show that the same conclusion holds for almost every $\Theta \in \mat$. 
 \end{prob}
 
 Notice that the weaker  inhomogeneous problem where the coprimality requirement $(\bq,\bp)\in P(\Pi) $ has been removed is in fact a well-known result.
 That is the inhomogeneous version of the Khintchine-Groshev theorem which follows from Theorem 1 of \cite{ScA} for $m=1$ 
 and Theorem 15 in Chapter 1 of \cite{Sp} when $m\ge 2$. We have established in Theorem 2 of \cite{DLN}
that Problem 1 holds true in the homogeneous case $\by =0$.
It is also shown in \cite{DLN} that Theorem 1 is equivalent to the following result concerning the smallness of  a generic
system of inhomogeneous linear forms:  

\begin{thm}
For every $\by\in\bR^n$ and for almost every  matrix 
$X \in {\rm Mat}_{n,m+n}(\bR) $, there exist infinitely 
many points $(\bq,\bp)\in P(\Pi) $
such that
$$
\left\vert X \left(\begin
    {matrix} 
     \bq\cr \bp
      \end{matrix}
      \right) -\by \right\vert \le \psi(| \bq |). $$
 \end{thm}
 The constant term $\by$ is fixed in Theorem 2, on the contrary to Theorem 1. 
 Finer statements of this kind, but without any coprimality constraint, have been recently
worked out by Dickinson,  Fischler, Hussain,  Kristensen and Levesley. We refer  to \cite{FHKL, HuKr}, and to the references therein,  for a discussion of their results.

\begin{prob}
Suppress, or relax, the (unnecessary ?) assumptions on $\Pi$ and $\psi$ occurring in the theorems 1 and 2.
 Namely,  the lower bounds $\card \, \pi_j\ge n+1$ and the monotonicity of the function $x\mapsto x^{m-1}\psi(x)^n$. 
 \end{prob}
 
 It is worth noting that the condition $\card \, \pi_j\ge n+1$ is needed to ensure the ergodicity of the  action of some group $\Gamma_\Pi$ used in the proof
 of Theorems 1 and 2.
  More precisely, let us denote by ${\rm SL}_{\pi_j}(\bZ)$ the subgroup of $\SLZ$ acting as an unimodular matrix on the coordinates with index in $\pi_j$ 
 and as the identity elsewhere.
 Let $\Gamma_\Pi= \prod_{j=1}^t {\rm SL}_{\pi_j}(\bZ)$ be the product  of these commuting subgroups.  Clearly, $\Gamma_\Pi$  is  a subgroup of $\SLZ$. 
 Then,  the action of $\Gamma_\Pi$ on $\Mat$ by right matrix multiplication is  ergodic, with respect to the invariant 
 Lebesgue measure on $\Mat \simeq \bR^{n(m+n)}$,  if and only if  $\card \, \pi_j\ge n+1$ for every  index $j= 1, \dots ,t$. 
 Notice in particular that  the action of $\SLZ$ on $\Mat$ is ergodic.

    \section{Exponent of Diophantine approximation of   dense orbits}
    
We  deal in this section with quantitative versions of the assertion (iv). To that purpose, we introduce the following exponent of Diophantine approximation.
  Let $X$ and $Y$ be two  matrices in $\Mat$. We define $\mu(X,Y)$ as the supremum of the real numbers $\mu$ such that there exist infinitely 
  many $\gamma \in \SLZ$ satisfying the inequality
  $$
  | X \gamma - Y | \le | \gamma |^{-\mu} .  \leqno{(1)}
  $$  
  Few informations are known concerning  the exponent $\mu(X,Y)$. In the case $m=n=1$, the following theorem follows from \cite{LN}, see also 
  \cite{MaWe} for a weaker lower bound in (i). We identify ${\rm Mat}_{1,2}(\bR)$ with  $\bR^2$.  
   \begin{thm}
   Let $X=(x_1,x_2)$ be a non-zero point in $\bR^2$ with irrational slope $x_1/x_2$.
   \\
   {\rm (i)} For every point $Y$ in $\bR^2$, we have the lower bound $\mu(X,Y) \ge 1/3$.
   \\
  {\rm(ii)}  If Y is any  non-zero point in $\bR^2$ with rational slope (i.e. the ratio of the two coordinates is a rational number), the equality $\mu(X,Y) =  1/2$ holds.
 \end{thm}
 
 Singhal \cite{Si} has partly extended the results of Theorem 3 to the linear action of the group ${\rm SL}_2(\cO_K)$ on $\bC^2$, where $\cO_K$ stands for the ring of integers of an imaginary quadratic field $K$ for which a convenient theory of continued fractions is available.
 
 Our knowledge concerning generic values of $\mu(X,Y)$ is a bit richer,  thanks to the recent results of \cite{GGNA,GGNB}. 
 Observe that the function $\mu(X,Y)$ is $\SLZ\times\SLZ$-invariant by componentwise right matrix multiplication. Since the action of $\SLZ$
 on $\Mat$ is ergodic, the function $\mu(X,Y)$  is equal to a constant $\mu_{m,n}$ (say) almost everywhere. In other words, for a fixed  exponent
 $\mu < \mu_{m,n}$ and for almost every pair $(X,Y)\in \Mat\times\Mat$,  the inequation (1) has infinitely many solutions $\gamma \in \SLZ$,
 while if $\mu  > \mu_{m,n}$  it   has only finitely many solutions $\gamma \in \SLZ$  almost surely. The determination of 
 the value of $\mu_{m,n}$ is an open problem. First, we  prove  the following upper bound:
 \begin{thm} 
 For any integers $m\ge 1$ and $n\ge 1$, we have the inequality
 $$
 \mu_{m,n} \le { m(m+n-1)\over n(m+n)}.
 $$
 \end{thm}
 
 \pro
 We follow the proof of the special case $m=n=1$ given in Section 5 of \cite{LNA}. For simplicity, put $a=m(m+n-1)$
 and $d =n(m+n)$. Let $\mu$ be a real number $> a/d$. We plan  to show that the set
 $$
 \cE_\mu = \left\{ (X,Y)\in \Mat\times \Mat \, ; \quad \mu(X,Y) > \mu \right\}
 $$
 has null Lebesgue measure.
  This will imply that $\mu_{m,n}\le \mu$ and next that $\mu_{m,n} \le a/d$ by letting $\mu$ tend to $a/d$.

 Denote by 
 $
 \Matmax 
 $ 
the open subset of $\Mat$ consisting of  the matrices with maximal rank $n$. The set $\Matmax$ is clearly stable by right multiplication by $\SLZ$.
Let $\Omega$ be a compact subset of $ \Matmax$ with smooth boundary and let $X \in \Matmax$ be such that the orbit 
 $X \SLZ$ is dense in $\Mat$. We claim that the set 
 $$
 \cE_\mu(X,\Omega) = \left\{ Y \in \Omega \, ; \quad (X,Y) \in \cE_\mu \right\}
 $$
is a null set. 
We naturally identify $\Mat$ with $ \bR^d$.
 For any point $Y\in \Mat$ and any positive real number $r$, let us  denote by
$$
B(Y,r) = \{ Z \in \Mat \, ; \quad | Z - Y | \le r\}
$$
the closed ball centered at $Y$ with radius $r$, whose Lebesgue measure $\lambda( B(Y,r))$ is equal to  $(2r)^d$.
 Now, let $Y$ belongs to $ \cE_\mu(X,\Omega)$. This means that there exist infinitely many $\gamma\in \SLZ$ satisfying (1). It follows that
$Y$ belongs to infinitely many balls $B(X\gamma,|\gamma|^{-\mu})$. Since $Y\in \Omega$, the center $X\gamma$ belongs to the neighborhood
$\Omega_{|\gamma|^{-\mu}}$ of points whose distance to $\Omega$ is $\le |\gamma|^{-\mu}$. Therefore, for any $\eta>0$, we have the inclusions
$$
 \cE_\mu(X,\Omega) \subseteq  \limsup_{{\gamma \in \SLZ\atop X\gamma \in \Omega_{\vert \gamma\vert^{-\mu}}}}B(X\gamma,| \gamma|^{-\mu})
\subseteq \limsup_{{\gamma \in \SLZ\atop X\gamma \in \Omega_\eta}}B(X\gamma,| \gamma|^{-\mu}) =: \cB_\eta.
 $$
Fix an $\eta>0$ sufficiently small so that $\Omega_\eta$ remains contained in the open set $\Matmax$.
We now use a fundamental counting result due to Gorodnick. Setting 
$$
N_R =  \Card \{ \gamma \in \SLZ \, ; \quad  X\gamma \in \Omega_\eta , | \gamma | =R \}, 
$$
  we have the upper bound
$$
N_1+ \cdots +N_R =  \Card \{ \gamma \in \SLZ \, ; \quad  X\gamma \in \Omega_\eta , | \gamma | \le R \} \le cR^{a}
\leqno{(2)}
$$
for some positive constant $c$ depending only on $X$ and $\Omega_\eta$ and any  positive integer $R$. 
Indeed, assuming that the orbit $X \SLZ$ is dense in $\Mat$, it is established in Theorem 3 of \cite{Go} that we have  an  asymptotic  equivalence of the form
$$
\Card \{ \gamma \in \SLZ \, ; \quad  X\gamma \in \Omega_\eta , | \gamma | \le R \}\sim  \delta R^{a}
$$
  as $R$ tends to infinity, with an explicit formula for the coefficient  $\delta$ depending only on $X$ and $\Omega_\eta$. 
  Using (2), we majorize
 $$
\begin{aligned}
 \sum_{{\gamma \in \SLZ\atop | \gamma | \le R \, ,\,  X\gamma \in \Omega_\eta}} & \lambda\Big(B(X\gamma,| \gamma|^{-\mu})\Big)
=
2^{d}\sum_{k= 1}^R {N_k \over k^{ d \mu }}
\\
  & =2^{d}\left( \sum_{k= 1}^{R-1} (N_1 + \cdots + N_k)\left( {1\over k^{d\mu}}- {1\over (k+1)^{d\mu}}\right) 
+ {N_1+ \cdots + N_R\over R^{d \mu}} \right)
\\
& \le 2^{d}c \left(  \sum_{k=1}^{R-1} k^{a} \left( {1\over k^{d\mu}} - {1\over (k+1)^{d\mu}}\right)  +  {  R^{a} \over R^{d\mu}}\right) 
\\
&\le 2^{d}c\mu d \left( \sum_{k=1}^{R-1} k^{a-d\mu -1}\right) + 2^{d}c R^{a-d\mu}.
\end{aligned}
$$
We deduce from the  above inequality that the sum 
$$
 \sum_{{\gamma \in \SLZ\atop   X\gamma \in \Omega_\eta}}  \lambda\Big(B(X\gamma,| \gamma|^{-\mu})\Big)
 $$
converges,  since we have assumed that $\mu > a/d$. By    Borel-Cantelli  lemma, the limsup set $\cB_\eta$ has null Lebesgue measure, as well as it subset 
$\cE_\mu(X,\Omega)$. We have proved the claim.

Selecting a countable family of  compacts $\Omega$ covering $ \Matmax$,  we deduce from the claim that the fiber
$
\{ Y \in \Mat ; (X, Y) \in \cE_\mu \}
$
of $\cE_\mu$ over $X$ has null Lebesgue measure whenever  the orbit $X \SLZ$ is dense. 
Now, Fubini theorem   yields that $\cE_\mu $ is a null set, since the  set of $X\in \Mat$ for which $X \SLZ$ is not dense
has  null Lebesgue measure   by the  criterion of Dani and Raghavan.
\cqfd

\medskip
\noindent
{\bf Remark: } The set  $\Matmax$ is an homogeneous space $H\backslash {\rm SL}_{m+n}(\bR)$, where $H$ is the semi-direct product
$H = {\rm SL}_{m}(\bR)\ltimes {\rm Mat}_{m,n}(\bR)$
and ${\rm SL}_{m}(\bR)$ acts on the $\bR$-vector space ${\rm Mat}_{m,n}(\bR)$ by left matrix multiplication. One can also prove Theorem 4 as part of  deep results
obtained by Ghosh, Gorodnick and Nevo in \cite{GGNA,GGNB}. 
They  give estimates of  generic exponents  of Diophantine approximation  associated to the action of lattice orbits on homogeneous spaces.
In fact, the above upper bound  corresponds  to the easiest  part of their results. See Section 3.1 of \cite{GGNB} and more specifically 
 Section 4.1 of \cite{GGNA} for details in the case $n=1$. 
 For comparison, notice however that the exponent $\kappa$ occurring in  \cite{GGNA,GGNB} is a {\it uniform} exponent,  according to the terminology of \cite{BuLa}, while our exponent $\mu$ is an ordinary  one. 
Following \cite{LN}, let us define the uniform variant  $\hat\mu(X,Y)$ of the exponent $\mu(X,Y)$  as the supremum of the real numbers
$\mu$ such that for any  large real number $R$, there exists $\gamma \in \SLZ$ satisfying 
$$
| \gamma | \le R \quad {\rm and} \quad | X \gamma - Y | \le R^{-\mu}.
$$
Obviously $\hat\mu(X,Y) \le \mu(X,Y)$ and it is expected that $\hat\mu(X,Y)$ and $\mu(X,Y)$
 have the same generic value $\mu_{m,n}$. The exponent $\kappa(X,Y)$  studied  in \cite{GGNA,GGNB}  is  the infimum of the real numbers
 $\kappa$ such the inequalities
$$
| \gamma | \le R^\kappa \quad {\rm and} \quad | X \gamma^{-1} - Y | \le R^{-1}.
$$
have a solution $\gamma \in \SLZ$ for any  large real number $R$. It is now expected that $\kappa(X,Y)$ has the generic  value $1/\mu_{m,n}$.
 
\begin{prob}
Show that the  formula 
$$
\mu_{m,n}= { m(m+n-1)\over n(m+n)} \leqno{(3)}
$$
 holds for every (at least one)  pair of positive integers $m$ and $n$.
\end{prob} 
 
In view of the  proof of Theorem 4, the conjecture may be considered as an optimistic version of the box principle.
 Fix a  compact subset $\Omega \subset \Matmax$.
 For large values of the radius  $R$, the set $\{ \gamma \in \SLZ \, ; \, X\gamma \in \Omega , | \gamma | \le R \} $ has cardinality $ \sim \delta R^{a}$
 whenever the orbit $X \SLZ$ is dense.
   It is expected that the images 
$X\gamma$ are well distributed in  $\Omega$ when $\gamma$ ranges over this set, at least  for a generic matrix  $X$. 

However,  Theorem 3 only yields the lower bound $\mu_{1,1} \ge 1/3$ and it follows from  Proposition 4.1 of \cite{GGNA} that $\mu_{m,1} \ge m/(m+1)$ when $m\ge 2$. 
Speculating further on the validity of formula (3), we shall prove in  Section 5 that $\mu_{1,1}$ has the expected value $1/2$, assuming  some strong properties
on the distribution of truncated ${\rm SL}_2(\bZ)$-orbits of primitive points in the plane, which remain unproved but seem plausible. This is the topic of the next section.

\section
 {Discrete orbits in the plane}
 
 The set of primitive points $P(\bZ^2)$ is an $\SL$-orbit in the plane. In other words,  we have $X\SL = P(\bZ^2)$ for any $X\in P(\bZ^2)$.
 We are interested in the repartition of the finite set of primitive points $X\gamma$ when $\gamma\in \SL$ ranges over  a ball $| \gamma | \le R$.
 Let us first recall a result due to Erd\H{o}s \cite{Er} on the distribution of  primitive points in the real plane: 
 \begin{Erdos}
 \noindent
 \\
 {\rm (i)} For any $Y$ in $\bR^2$ with $| Y|$ large enough, there exists $Z \in P(\bZ^2)$ such that
 $$
 | Z - Y | \le  { \log | Y | \over \log\log | Y|}.
 $$
 \\
{\rm (ii)} There exists  $Y$ in $\bR^2$ with $| Y|$ arbitrarily large such that  we have the lower bound
 $$
 | Z - Y | \ge {1\over 2} \left({ \log | Y | \over \log\log | Y|}\right)^{1/2},
 $$
for any  $Z \in P(\bZ^2)$.
 \end{Erdos}

  Let $Y =(y_1,y_2)$ be the coordinates of $Y$. Erd\H{o}s' results are  in fact expressed differently in term of the quantity $\min( |y_1|,|y_2|)$ rather than $| Y|$.  
  The above formulation  involving the norm $| Y|= \max( |y_1|,|y_2|)$ is a straightforward  consequence of the statements (2) and (3) of \cite{Er}
  and their proof. 
 An interesting, but difficult,  open question  is to sharpen  (i) and/or  (ii):
 
 \begin{prob}
 Determine the smallest possible upper bound in the assertion {\rm (i)} of  Erd\H{o}s Theorem as a function of the norm $| Y |$. 
 \end{prob}
  
 We now come to the study of the orbit $X \SL$ where $X$ is a given primitive point. 
For each point $Z \in P(\bZ^2)$, we choose  a matrix ${\bf Z}\in \SL$ whose first row is equal to $Z$ having norm 
  $|{\bf Z}|=|Z|$. Using B\'ezout, it is easily seen that there are two such matrices ${\bf Z}$, unless $Z=(\pm 1,0)$ or $Z= (0,\pm 1)$  in which cases
  there are three.   Then, we can obviously reformulate  (i) in the matricial form
  $$
   |(1,0) {\bf Z} - Y | \le  { \log | Y | \over \log\log | Y|}.
 $$
 Writing $ (1,0) = X \bX^{-1}$, we replace in the above inequality the base point $(1,0)$ by an arbitrary primitive point $X\in P(\bZ^2)$. 
 Setting  $\gamma = {\bf X}^{-1}{\bf Z}\in \SL$ and $ R = 3 | X| |Y|$, we  thus obtain the estimates
 $$
 | \gamma| \le 2 | X | \left(| Y| +   {  \log | Y | \over \log\log | Y|}\right) \le R \quad {\rm and} \quad | X \gamma -Y | \le   { \log  R  \over \log\log  R},
 \leqno{(4)}
$$ 
when $| Y |$ is large enough.
We wish to find analogous estimates for smaller values of $R$,
enlarging possibly the  upper bound $  \log  R  / \log\log  R $ for  the distance between $X\gamma$ and $Y$.
More precisely we propose the following

 \begin{prob}
  Let $\varepsilon$ be a positive real number.
 Show that for  any large real number $R$, any primitive point
 $X \in P(\bZ^2)$ and any real point  $Y\in \bR^2$
 with norm $| Y | \le (| X | R)^{{1\over 1+\varepsilon}}$, there exists $\gamma\in \SL$ such that 
 $$
 |\gamma| \le R \quad {\rm and} \quad | X \gamma -Y | \le | X |^{1+\varepsilon} R^{\varepsilon}.\leqno{(5)}
 $$
  \end{prob}

The above hypothesis  may be viewed as a property of uniform  repartition of the truncated orbit $\{ X\gamma \, ; \, \gamma\in \SL, | \gamma | \le R\}$ in the following way.
The number of $\gamma \in \SL$ whose norm $| \gamma | $ is $\le R$ is  $\asymp R^2$, while $| X \gamma | \le 2 | X | | \gamma | \le 2| X | R$.
If we divide the square $\{ Y \in \bR^2; | Y | \le  2 | X | R\}$ into $\asymp R^2$ small squares with side $\asymp | X |$, each small square should ideally
contain a point  of the set $\{ X\gamma \, ; \, \gamma\in \SL, | \gamma | \le R\}$.
On the other hand, the correspondence $\gamma\mapsto X\gamma$ is not one to one, but the objection is irrelevant here.
Indeed, putting  $Z =X\gamma \in P(\bZ^2)$,  we can express
$\gamma \in \SL$ in the form
$$
\gamma =
\left(\begin{matrix}
x'_2z_1 - x_2z'_1-u x_2 z_1& x'_2z_2 -x_2z'_2 - ux_2z_2
\\
-x'_1z_1+x_1z'_1+u x_1z_1 & -x'_1z_2+x_1z'_2+ux_1z_2
\end{matrix}
\right)
$$
for some $u\in \bZ$, where
$$
{\bf X}= \left(\begin{matrix}
x_1 &x_2 \\
x'_1 &x'_2 
\end{matrix}
\right)
\quad{\rm and} \quad
{\bf Z}= \left(\begin{matrix}
z_1 &z_2 \\
z'_1 &z'_2 
\end{matrix}
\right)
$$
are the respective liftings of $X$ and $Z$ in $\SL$ previously considered. It follows that there exist at most $2R/(|X| |Z|)+1$
matrices $\gamma \in \SL$ for which $|\gamma | \le R$ and $X\gamma =Z$. The number of such matrices $\gamma$ is thus bounded by  3 when 
$R \le  |X| |Z|$. Now, the number of $\gamma\in \SL$ for which  $| \gamma | \le R$ and $| X \gamma | \le R /| X |$ is at most
$$
\sum_{ { Z \in P(\bZ^2)\atop | Z | \le R / | X|}}   {2 R\over | X | | Z|} +1 
\le  {3 R\over | X|}\sum_{ { Z \in P(\bZ^2)\atop | Z | \le R / | X|}} {1\over | Z|}
= {3 R\over | X|} \sum_{k=1}^{\lfloor R/| X |\rfloor} {8 \varphi(k)\over k} \le {24R^2 \over | X |^2 }.
$$
We easily deduce from the preceding considerations that
 the cardinality of the set of points $\{ X\gamma \, ; \, \gamma\in \SL, | \gamma | \le R\}$ is $\asymp R^2$ when $R$ is sufficiently large independently of $X$,
 where the two multiplicative coefficients implicitely  involved in the symbol $\asymp$ are absolute constants.

Notice that   we have relaxed somehow the constraints by 
assuming that 
$$
| Y |  \le (| X | R)^{1/(1+\varepsilon)}\le 2 | X | R
$$
 and requiring only the weaker upper bound 
$| X \gamma -Y | \le | X |^{1+\varepsilon} R^{\varepsilon}$,
  because  these  estimates  are sufficient for our purpose. Of course, sharper estimates may reveal useful as well.

\smallskip
\noindent
{\bf Remark.}
 The finite set  $\{ X\gamma \, ; \, \gamma\in \SL, | \gamma | \le R\}$
does not look like a grid of step $| X|$ in $\bR^2$, as one might believe at first glance. 
Any point $Y\in P(\bZ^2)$ with norm $| Y | \le R/(2 | X|)$
belongs to this  set (take $\gamma= \bX^{-1}\bY$ as in (4)). Therefore, we have a concentration of points around the origin when $R \gg | X|$.
However,  the analogy with a grid remains meaningful inside  some annulus  centered at the origin.

 \section{A corollary}
 
 We relate  the distribution of truncated integral orbits to the rate of density  of dense orbits, 
  showing   that a positive answer to Problem 5 implies  that $\mu_{1,1}= 1/2$.
   
   We say that a real number   $y$ is  {\it  very well approximated by rationals} 
  if for some  exponent $\omega >2$, the inequation $| y-p/q| \le q^{-\omega}$ has infinitely many rational solutions $p/q$. 
 An irrational number $y$ is not very well approximated by rationals if and only if the sequence $(s_j)_{j\ge 0}$
 of denominators of its convergents satisfies the asymptotic growth condition
 $s_{j+1}\le s_j^{1+\varepsilon}$ for every  $\varepsilon>0$. It is well-known that  almost every real number $y$ is not very well approximated by rationals.
  It then suffices to prove the 
 
 \begin{prop}
 Let  $X=(x_1,x_2)$  and  $Y= (y_1,y_2)$ be two points in $\bR^2$.  We suppose that  the point $Y$ does not belong to the orbit $X \SL$, that the slope $ x_1/x_2$ of $X$ is an irrational number and that the slope $ y_1/y_2$ of $Y$ is an irrational number which  is  not  very well approximated by rationals. 
  Assume  that Problem 5 has been affirmatively resolved  for any 
 $\varepsilon >0$. Then $\mu(X,Y)\ge 1/2$.
 \end{prop}

\pro 
We have to show that
 for any $\eta >0$, there exist infinitely many  $\gamma \in \SL$ such that
 $$
 | X \gamma  - Y | \le | \gamma |^{-1/2 + \eta}.\leqno{(6)}
 $$
The proof  is a variant  of that of Theorem 1 in \cite{LN}: we replace the unipotent matrix of the form $ \left(\begin{matrix}
1 &\ell\\
0 &1 
\end{matrix}
\right)
$ occurring there by a matrix $G$ whose transposed matrix  is a solution to Problem 5.  

Put $\xi = x_1/x_2$ and $y=y_1/y_2$. Let   $(p_k/q_k)_{ k\ge 0}$ and $(t_j/s_j)_{ j\ge 0}$ 
   be  the sequences of  convergents of $\xi$ and $y$ respectively and set  
$$
  M_k = \left(\begin{matrix} q_k & q_{k-1}\\ -p_{k}& -p_{k-1}\end{matrix}\right),
  \quad 
  N_j = \left(\begin{matrix} t_j & s_j \\ t_{j-1} & s_{j-1}\end{matrix}\right).
 $$
The indices $k$  and $j$ will be restricted to  odd integers, so that both $M_k$ and  $N_j$  belong to $\SL$.
We construct $\gamma $ in  the form
 $$
 \gamma = M_k G N_j
 $$
 where 
 $$
   G =  \left(\begin{matrix} a & b \\ c & d\end{matrix}\right) \in \SL
 \quad {\rm and } \quad 
 | G | \le R.
 $$
 The indices $j,k$ and the norm upper bound $R$ will be chosen later. 
 
 Put $\epsilon_n = q_n \xi -p_n$  and recall the estimate 
 $$
 {1\over 2 q_{n+1}} \le | \epsilon_n | \le {1\over q_{n+1}}
 $$
  derived  from the theory of continued fractions.
 Our starting point is the identity
 $
   q_k\epsilon_{k-1} - q_{k-1}\epsilon_k=1.
 $
 Multiplying by  $y_2/x_2$, we find 
 $$
 {y_2 \over x_2}=   (q_k {y_2\over x_2})\epsilon_{k-1} + (-q_{k-1} {y_2\over x_2}) \epsilon_k \leqno{(7)}
 $$
 and on the other hand, we have the obvious equality
 $$
 0 =  (q_k {y_2\over x_2})q_{k-1} + (-q_{k-1} {y_2\over x_2}) q_k. \leqno{(8)}
 $$
 Now write
 $$
 X\gamma -Y = x_2(\xi,1) M_kGN_j -Y = x_2(\epsilon_k, \epsilon_{k-1})G N_j -Y =(z_1,z_2)
 $$
 with
 $$
 \begin{aligned}
 z_1 =&  x_2\Big((a t_j  + b t_{j-1} )\epsilon_k + (c t_j  +d t_{j-1})\epsilon_{k-1}\Big) - y_1
 \\
 z_2 =& x_2\Big((a s_j +b s_{j-1} )\epsilon_k  + (c s_j  +d s_{j-1}Ê )\epsilon_{k-1}\Big) -y_2
 \\
 =& x_2\left(\left(  a s_j +b s_{j-1} + q_{k-1}{y_2\over x_2}\right)\epsilon_k + \left(  c s_j +d s_{j-1} - q_k{y_2\over x_2}\right)\epsilon_{k-1}\right),
  \end{aligned}
  \leqno{(9)}
$$ 
using (7) for the last equality. Put 
$$
 \Delta = \max \left(\left\vert a s_j +b s_{j-1} + q_{k-1}{y_2\over x_2}  \right\vert , \left\vert c s_j +d s_{j-1} - q_k{y_2\over x_2 } \right\vert \right).
 $$
We  immediately  deduce  from  (9) that 
$$
| z_2 | \le 2 | x_2|  {\Delta\over q_k}.\leqno{(10)}
$$
 Using again the expressions (9), observe now that 
$$
  | z_1 - y z_2 | 
  =   | x_2| \Big\vert(t_j-s_j y)(a\epsilon_k + c \epsilon_{k-1}) + (t_{j-1}-s_{j-1} y)(b \epsilon_k +d  \epsilon_{k-1})\Big\vert
 \le 4 | x_2| { R\over s_j q_k}.\leqno{(11)}
 $$
Combining (10) and (11), we deduce from the triangle inequality  the upper bound
$$
| X \gamma - Y | = \max (| z_1|, | z_2|) \ll {\Delta\over q_k} + { R\over s_j q_k}.
\leqno{(12)}
$$

We now bound the norm of $\gamma= M_kGN_j$ and  claim that
$$
| \gamma | \ll {\Delta q_k}+{s_jR\over q_k}+  {q_kR \over s_j}.
\leqno{(13)}
$$
Expand the product
$$
 \gamma = \left(\begin{matrix} q_k & q_{k-1}\\ -p_{k}& -p_{k-1}\end{matrix}\right)
  \left(\begin{matrix} a & b \\ c & d\end{matrix}\right)
 \left(\begin{matrix} t_j & s_j \\ t_{j-1} & s_{j-1}\end{matrix}\right)
 =\left(\begin{matrix} A & B\\ C& D\end{matrix}\right),
  $$
  where we have set
  $$
  \begin{aligned}
  A  & = q_k(at_j+ bt_{j-1}) + q_{k-1}(ct_j+dt_{j-1}),
  \\B  &= q_k(as_j+ bs_{j-1}) + q_{k-1}(cs_j+ds_{j-1}),
   \\
C & = -p_k(at_j+ bt_{j-1}) - p_{k-1}(ct_j+dt_{j-1}),
\\
  D &  = -p_k(as_j+ bs_{j-1}) - p_{k-1}(cs_j+ds_{j-1}).
\end{aligned}
$$
Since 
$$
(t_j,t_{j-1}) = y(s_j, s_{j-1}) + \cO\left({1\over s_j}\right) \quad {\rm and } \quad( p_k,p_{k-1}) = \xi(q_k, q_{k-1}) + \cO\left({1\over q_k}\right),
$$
we have the estimates
$$
A= y B + \cO\left({Rq_k\over s_j}\right), D = -\xi B + \cO\left({Rs_j\over q_k}\right), C =-\xi y B +\cO\left({Rq_k\over s_j}\right)+\cO\left({Rs_j\over q_k}\right).
$$
Using  (8), we can express $B$ in the form
$$
B= q_k\left(as_j+ bs_{j-1} + q_{k-1}{y_2\over x_2}\right) + q_{k-1}\left(cs_j+ds_{j-1}- q_{k}{y_2\over x_2}\right)
$$
which yields the bound $| B | \le 2 q_k \Delta$.  The inequality  (13) is thus established.

At this stage, we use the  (conjectural) estimate (5) of the preceding section applied to the pair of points
$$
  (s_j,s_{j-1})\in P(\bZ^2) \quad{\rm and}\quad  (-q_{k-1}{y_2\over x_2}, q_{k}{y_2\over x_2})\in \bR^2
 $$
 in order to majorize  $\Delta$ in a non-trivial way. Fix arbitrarily $\varepsilon >0$ and set 
$$
R= { (| y_2| q_k/| x_2|)^{1+ \varepsilon} \over  s_j }.\leqno{(14)}
$$
The condition
 $$
  { q_k | y_2 | \over | x_2 | }= \left\vert  (-q_{k-1}{y_2\over x_2}, q_{k}{y_2\over x_2})\right\vert \le \left( | (s_j,s_{j-1})| R\right)^{{1\over 1+ \varepsilon}}
  =  \left( s_jR\right)^{{1\over 1+ \varepsilon}}
  $$
  occurring in Problem 5 is obviously satisfied.
 It then follows from (5) that there exists 
 a matrix  $G =  \left(\begin{matrix} a & b \\ c & d\end{matrix}\right) \in \SL$ with norm
$| G | \le R$ such that
$$
 \Delta  = \max \left(\left\vert a s_j +b s_{j-1} + q_{k-1}{y_2\over x_2}  \right\vert , \left\vert c s_j +d s_{j-1} - q_k{y_2\over x_2 } \right\vert \right) \le s_j^{1+\varepsilon} R^{\varepsilon}.\leqno{(15)}
 $$
Now fix an odd   index $k$ and choose $j$ odd so that $s_j$ should be located in the interval
$$
q_k^{1/3} \le s_j \le q_k^{1/3 +\varepsilon}.\leqno{(16)}
$$
This is possible when  $k$ is large enough,  since we have assumed that $y$ is not very well approximated by rationals. Then, 
 combining the estimates (12) to (16) easily yields the upper bounds
$$
| X \gamma -Y  | \ll q_k^{-2/3 + 2 \varepsilon + 2 \varepsilon^2}
 \quad\hbox{\rm and  } \quad  | \gamma | \ll q_k^{4/3 + 2\varepsilon + 2\varepsilon^2}.
 $$
 It follows that 
 $$
 | X \gamma -Y  |\ll q_k^{-2/3 +2 \varepsilon + 2 \varepsilon^2 } \ll | \gamma |^{-1/2 + \eta},
 $$
 where  $\eta = 9(\varepsilon+ \varepsilon^2)/(4+6\varepsilon+ 6\varepsilon^2)$ is arbitrarily small when $\varepsilon$ is sufficiently small.
  We have thus established (6).
  
  The above construction produces infinitely many solutions $\gamma$ to (6) when $k$ ranges over the odd integers, because 
 the distance  $| X \gamma -Y  |$ tends  to $0$ as $k$ tends to infinity,  while remaining  positive 
 (recall that we have assumed that $Y$ does not belong to the orbit $X\SL$). The proposition is proved.\cqfd

 {\small

 Institut de Math\'ematiques de Marseille, Aix Marseille Universit\'e,  Case 907, 163 avenue de Luminy, 13288,  Marseille C\'edex 9.

  {\tt michel-julien.laurent@univ-amu.fr}, 
  }

\end{document}